\documentclass[letterpaper, 10 pt, conference]{ieeeconf}  

\IEEEoverridecommandlockouts                              

\usepackage[top=1in, bottom=1in, left=0.75in, right=0.75in]{geometry}
\usepackage{amsmath, amssymb, epsfig, graphicx, bm}
\usepackage{cases}
\usepackage{float, subfigure}
\usepackage{color}
\usepackage{amsthm}
\usepackage{amsfonts}
\usepackage{cite, url}
\usepackage[all,cmtip]{xy}
\usepackage{color,hyperref}
\definecolor{darkblue}{rgb}{0,0,1}
\hypersetup{colorlinks,breaklinks,
linkcolor=darkblue,urlcolor=darkblue,anchorcolor=darkblue,citecolor=darkblue}

\usepackage[normalem]{ulem} 

\newcommand{\pro}[2]{{\mathcal{P}_{#1}\left\{{#2}\right\}}}

\newcommand{\tanya}[1]{\normalsize{{\color{black}\ #1}}}

\newtheorem{theorem}{Theorem}
\newtheorem{definition}{Definition}
\newtheorem{proposition}{Proposition}
\newtheorem{lemma}{Lemma}
\newtheorem{remark}{Remark}

\newtheorem{assumption}{Assumption}
\newtheorem{corollary}{Corollary}

\newcommand{\R}{{\mathbb{R}}}
\newcommand{\bx}{{\mathbf{x}}}
\newcommand{\hbx}{\hat{\bx}}

\newcommand{\di}{{\mathrm{diag}}}

\newcommand{\dW}{{\|I-W\|^2}}

\newcommand{\dJi}[1]{{\nabla_{#1} J_{#1}}}

\newcommand{\A}{{\mathcal{A}}}

\newcommand{\Gra}{{\mathcal{G}}}

\newcommand{\tx}{{\tilde{x}}}

\newcommand{\Om}{{\Omega}}

\newcommand{\one}{{\mathbf{1}}}

\newcommand{\N}{{\mathcal{N}}}

\newcommand{\bF}{{\mathbf{F}}}

\newcommand{\trace}{{\mathrm{trace}}}
\newcommand{\T}{{\mathrm{T}}}

\def\an#1{{\color{black}#1}}

\title{\LARGE \bf
Accelerating Distributed Nash Equilibrium Seeking}

\author{Tatiana Tatarenko$^{1}$ and Angelia Nedi\'c$^{2}$
\thanks{*This work was not supported by any organization}
\thanks{$^{1}$Tatiana Tatarenko is with the Department of Control Theory and Intelligent Systems, TU Darmstadt, German
        {\tt\small tatiana.tatarenko@tu-darmstadt.de}}%
\thanks{$^{2}$Angelia Nedi\'c is  with School of Electrical, Computer and Energy Engineering,
        Arizona State University, USA.
        {\tt\small Angelia.Nedich@asu.edu}. Te work has been done under support of the Office of Naval Research award N00014-21-1-2242.%
}
}

\begin{document}

\maketitle
\thispagestyle{empty}
\pagestyle{empty}

\begin{abstract}
This work proposes a novel distributed approach for \an{computing a Nash equilibrium} in convex games with restricted strongly monotone pseudo-gradients. By leveraging the idea of the centralized operator extrapolation method presented in \cite{LanExtrapolation} to solve variational inequalities, we develop the algorithm converging to Nash equilibria in games, where players have no access to the full information but are able to communicate  with neighbors over some communication graph. The convergence rate is demonstrated to be geometric and improves the rates obtained by the previously presented procedures seeking Nash equilibria in the class of games under consideration.

\end{abstract}

\section{INTRODUCTION}\label{sec:intro}

Game theory deals with a specific class of optimization problems arising in multiagent systems, in which each
agent, also called player, aims to minimize its local cost function coupled through decision variables (actions) of all agents (players) in a system. 
The applications of game-theoretic optimization can be found, for example, in electricity markets, communication networks, autonomous driving systems and the future smart grids \cite{Alpcan2005, BasharSG, GTM, Scutaricdma}. Solutions to such optimization problems are Nash equilibria  which characterize desirable and stable joint actions in games. To find these solutions in \an{a so called convex game, one can use their equivalent characterization as the solutions to the variational inequality} defined for the game's pseudo-gradient over the joint action set \cite{FaccPang1}. Moreover, it is known that, given a strongly monotone and Lipschitz continuous \an{mapping, the projection algorithm} converges geometrically fast to the unique solution of the variational inequality  and, thus, to the unique Nash equilibrium \an{of} the game. 
\an{The convergence rate, in terms of the $k$th iterate's distance to the solution,
is in the order of $O\left(\exp\left\{-\frac{k}{\gamma^2}\right\}\right)$(see \cite{Nesterov}), where 
$\gamma = L/\mu\ge 1$ with $L$ and $\mu$ being the Lipschitz continuity and strong monotonicity constants of the mapping, respectively.}
\an{This rate has been improved in \cite{Nesterov} 
to the rate of $O\left(\exp\left\{-\frac{k}{\gamma}\right\}\right)$ 
by a more sophisticated algorithm
that requires, at each iteration, two operator evaluations and two projections}. To relax these requirements, the paper \cite{LanExtrapolation} presents the so called operator extrapolation method achieving the same rate $O\left(\exp\left\{-\frac{k}{\gamma}\right\}\right)$ with one
operator evaluation and one projection per iteration. Moreover, geometrically fast convergence of the operator extrapolation method takes place under a weaker condition of restricted strong
monotonicity (see \textbf{Notations}).  
\an{However, these fast algorithms require full information in the sense that each player observes actions} of all other players at every iteration.

Since in the modern large-scale
systems each agent has access only to some partial information about joint actions, \emph{fast distributed communication-based} optimization procedures in games have
gained a lot of attention over the recent years (see \cite{survey_distGT} for an extensive
review and bibliography). In particular, the work \cite{Bianchi2019} presents a proximal-point algorithm for converging to the Nash equilibrium with
a geometric rate. However, this algorithm requires \an{the evaluation of a proximal operator,} at each iteration, and cannot be rewritten as iterations that give the next state in terms of the current one. On the other hand, the papers \cite{ifac_TatNed} and \cite{directMethod_Grammatico} propose the distributed procedures based on the gradient algorithm and demonstrate their geometric convergence \an{rate in the order $O\left(\exp\left\{-\frac{k}{\gamma^4}\right\}\right)$ for strongly monotone games with player communications
over time-invariant and time-varying graphs,} respectively. The works \cite{Cdc2018_TatShiNed, AccGRANE_TAC} focus on a reformulation of \an{ a Nash equilibrium in distributed setting in terms of a so called augmented variational inequality, which takes into account the communication network that players are using}. The main goal of such reformulations has been to adjust the fast centralized procedure from \cite{Nesterov} to the distributed settings and accelerate learning Nash equilibria in distributed settings. However, the acceleration (to the rate $O\left(\exp\left\{-\frac{k}{\gamma^3}\right\}\right)$) has been guaranteed only for a restrictive subclass of games with strongly monotone and Lipschitz continuous pseudo-gradients. This restriction is due to the fact that the mapping defining the augmented variational inequality is generally not strongly monotone but \an{it can be made} restricted strongly monotone. Moreover, the augmented variational inequality requires introduction of an extra parameter which has to be properly set up to guarantee convergence of the proposed algorithms.  

This article presents a novel fast distributed discrete-time algorithm for seeking Nash equilibria in games with restricted strongly monotone pseudo-gradients. To avoid issues related to the augmented variational inequality arising in the distributed settings and still to be able to accelerate the previously known rates, this algorithm leverages the idea of the operator extrapolation method from \cite{LanExtrapolation} instead of the Nesterov's acceleration approach presented in \cite{Nesterov}. 
We develop a procedure converging to the Nash equilibrium with the rate $O\left(\exp\left\{-\frac{k}{\gamma^2}\right\}\right)$ and requiring one projection and gradient calculation per iteration.

\textbf{Notations.}
The set $\{1,\ldots,n\}$ is denoted by $[n]$.
For any function $f:K\to\R$, $K\subseteq\R^n$, $\nabla_i f(x) = \frac{\partial f(x)}{\partial x_i}$ is the partial derivative taken in respect to the $i$th coordinate of the vector variable $x\in\R^n$.
For any real vector space $\tilde E$ its dual space is denoted by $\tilde E^*$ and the inner product is denoted by $\langle u,v \rangle$, $u\in\tilde E^*$, $v\in \tilde E$.
\an{A mapping $g:\tilde E\to \tilde E^*$ is said to be \emph{strongly monotone with the constant $\mu>0$ on the set $Q\subseteq \tilde E$}, if $\langle g(u)-g(v), u - v \rangle\ge\mu\|u - v\|^2$ for all $u,v\in Q$.
} 
\an{It is said to be
{\it restricted strongly monotone with respect to $u^*\in Q$}, if $\langle g(u)-g(u^*), u - u^* \rangle\ge\mu\|u - u^*\|^2$ for all $u\in Q$.}
We consider real vector space $E$, which is either space of real vectors $E = E^* = \R^n$ or the space of real matrices $E = E^* = \R^{n\times n}$. 
\an{In the case $E = \R^n$ we use $\|\cdot\|$ to denote the Euclidean norm induced by the standard dot product in $\R^n$. In the case $E = \R^{n\times n}$, 
the inner product $\langle u,v \rangle \triangleq \sqrt{\trace(u^Tv)}$ is the Frobenius inner product on $\R^{n\times n}$ and 
$\|\cdot\|$ denotes the Frobenius norm induced by the Frobenius inner product, i.e., $\|v\| \triangleq \sqrt{\trace(v^Tv)}$.}
We use $\pro{\Om}{v}$ to denote the projection of $v\in E$ on a set $\Om\subseteq E$.
For any matrix $A$, \an{the vector of diagonal entries of the matrix $A$ is denoted by $\di(A)$.}

\section{Distributed Learning in Convex Games}

We consider a non-cooperative game between $n$ players. Let $J_i$ and $\Om_i\subseteq \R$ denote\footnote{All results below are applicable for games with different dimensions $\{d_i\}$ of the action sets $\{\Om_i\}$. The one-dimensional case is considered for the sake of notation simplicity.} respectively the cost function and the feasible action set of the player $i$. We denote the joint action set by $\Om = \Om_1\times\cdots\times\Om_n$. Each function $J_i(x_i,x_{-i})$, $i\in[n]$, depends on $x_i$ and $x_{-i}$, where $x_i\in\Om_i$ is the action of the player $i$ and $x_{-i}\in\Om_{-i}=\Om_1\times\cdots\times\Om_{i-1}\times\Om_{i+1}\times\cdots\times\Om_n$ denotes the joint action of all players except for the player $i$. We assume that the players can interact over an undirected communication graph $\Gra([n],\A)$. The set of nodes is the set $[n]$ of players, and 
the set $\A$ of undirected arcs is such that $\{i,j\}\in\A$ whenever there is an undirected communication link between $i$ to $j$
and, thus, some information (message) can be passed between the players $i$ and $j$.
For each player $i$, the set $\N_i$ is the set of neighbors in the graph $\Gra([n],\A)$, i.e.,
\an{$\N_{i}\triangleq\{j\in[n]: \, \{i,j\}\in\A\}$.}
We denote this game by $\Gamma(n,\{J_i\},\{\Om_i\},\Gra)$, and 
we make the following assumptions regarding the game.

\begin{assumption}\label{assum:convex}[Convex Game]
For all $i\in[n]$, the set $\Om_i$ is convex and closed, while the function $J_i(x_i, x_{-i})$ is convex and continuously differentiable in $x_i$ for each fixed $x_{-i}$.
\end{assumption}
When the cost functions $J_i(\cdot,x_{-i})$ are differentiable, we can define the pseudo-gradient.
\begin{definition}\label{def:gamemapping} The \emph{pseudo-gradient} $F(x):\Om\to\R^n$ of the game $\Gamma(n,\{J_i\},\{\Om_i\},\Gra)$ is defined as follows:
 $F(x)\triangleq\left[\nabla_1 J_1(x_1,x_{-1}), \ldots, \nabla_n J_n(x_n,x_{-n})\right]^T\in\R^{n}$,
 \an{where $\nabla _i$ denotes the partial derivative with respect to $x_i$ (see \bf{Notations}).}
\end{definition}


A solution to a game  is a Nash equilibrium, defined below.
\begin{definition}\label{def:NE}
 A vector $x^*=[x_1^*,x_2^*,\cdots, x_n^*]^T\in\Om$ is a \emph{Nash equilibrium} if for all $i\in[n]$ and all $x_i\in \Om_i$
 $$J_i(x_i^*,x_{-i}^*)\le J_i(x_{i},x_{-i}^*).$$
 \end{definition}

By Assumption~\ref{assum:convex} and the connection between 
Nash equilibria and solutions of variational inequalities~\cite{FaccPang1}, 
the point $x^*\in\Om$ is a Nash equilibrium of the game $\Gamma(n,\{J_i\},\{\Om_i\},\Gra)$ if and only if the following
variational inequality holds
\begin{align}\label{eq:NE}
 \langle F(x^*), x-x^*\rangle\ge 0 \quad \mbox{for all $x\in\Om$}.
\end{align}
We make further assumptions regarding the players cost functions, as follows.
\begin{assumption}\label{assum:str_monotone}
 The game $\Gamma(n,\{J_i\},\{\Om_i\},\Gra)$ has a Nash equilibrium $x^*$, and the pseudo-gradient mapping $F(x)$ is defined on the whole space $\R^n$ and is \emph{restricted strongly monotone with respect to $x^*$} on $\R^n$ with a constant $\mu>0$.
\end{assumption}
The existence of a Nash equilibrium is guaranteed if, for example, Assumption~\ref{assum:convex} holds and the action sets $\Om_i$, $i\in[n]$, are bounded \cite{FaccPang1}. 
\begin{assumption}\label{assum:Lipschitz}
For every $i\in[n]$ the function $\dJi{i}(x_i,x_{-i})$ is Lipschitz continuous in
$x_i$ on $\Om_i$ for every fixed $x_{-i}\in\R^{n-1}$, that is, 
there exist a constant $L_i\ge 0$ such that  for all $x_{-i}\in\R^{n-1}$ we have for all $x_i,y_i\in\Om_i$,
\begin{align*}
 |\dJi{i}(x_i,x_{-i})-\dJi{i}(y_i,x_{-i})|&\leq L_i|x_i-y_i|.
\end{align*}
Moreover, for every $i\in[n]$ the function $\dJi{i}(x_i,x_{-i})$ is Lipschitz continuous in $x_{-i}$ on $\R^{n-1}$, for every fixed $x_i\in\Om_i$, that is, there is a constant $L_{-i}\ge 0$ such that 
for all $x_{i}\in\Om_i$ we have for all $x_{-i},y_{-i}\in \R^{n-1},$
\begin{align*}
|\dJi{i}(x_i,x_{-i})-\dJi{i}(x_i,y_{-i})|&\leq L_{-i}\|x_{-i}-y_{-i}\|.
\end{align*}
\end{assumption}

\an{The players' communications are restricted to the underlying connectivity graph $\Gra([n],\A)$, with which we associate a nonnegative symmetric mixing matrix $W$, i.e.,  a symmetric matrix with nonnegative entries and with positive entries $w_{ij}$ only when $\{i,j\}\in\A$. To ensure sufficient information  "mixing" in the network, we assume that the graph is connected. These assumptions are formalized, as follows.} 

\begin{assumption}\label{assum:connected}
The underlying undirected communication graph $\Gra([n],\A)$ is connected. The associated non-negative \an{symmetric} mixing matrix $W=[w_{ij}]\in\R^{n\times n}$ defines the weights on the undirected arcs such that $w_{ij}>0$  
if and only if $\{i,j\}\in\A$ and $\sum_{j=1}^{n}w_{ij} = 1$ for all $i\in[n]$.
\end{assumption}

\begin{remark}
There are some simple strategies for generating symmetric mixing matrices over undirected graphs for which Assumption~\ref{assum:connected} holds (see Section 2.4 in~\cite{Shi2014} for a summary of such strategies).
\end{remark}
Assumption~\ref{assum:connected} implies that the second largest singular 
value $\sigma$ of $W$ is such that $\sigma\in(0, 1)$ and for any $x\in\R^n$ the following average property holds (see \cite{OlshTsits}):
\begin{align}\label{eq:sigma}
	\|Wx-\one\bar{x}\|\le \sigma\|x-\one\bar{x}\|,
\end{align}
where $\bar{x} = \frac{1}{n}\langle \one,x\rangle$ is the average of the coordinates of $x$.

In this work, we are interested in \emph{distributed seeking of the Nash equilibrium} in a game $\Gamma(n,\{J_i\},\{\Om_i\},\Gra)$ for which Assumptions~\ref{assum:convex}--\ref{assum:connected} hold. We note that under Assumptions~\ref{assum:convex}--\ref{assum:str_monotone}, the game has a unique Nash equilibrium (see~\cite{AccGRANE_TAC}). 

\section{Algorithm Development}\label{sec:algo}

\subsection{Direct Acceleration}
Throughout the paper, we let player $i$ hold a \emph{local copy} of
the global decision variable\footnote{Note that global decision variable $x$ is a fictitious variable which never exists in the designed decentralized computing system.} $x$, which is denoted by
$$x_{(i)}=[\tx_{(i)1};\ldots;\tx_{(i)i-1};x_i;\tx_{(i)i+1};\ldots;\tx_{(i)n}]\in\R^n.$$
Here $\tx_{(i)j}$ can be viewed as a temporary estimate of $x_j$ by player $i$. \an{In this notation,} we always have $\tx_{(i)i}=x_i$. Also, we compactly denote the temporary estimates that player $i$ 
has for all decisions of the other players as
$$\tx_{-i}=[\tx_{(i)1};\ldots;\tx_{(i)i-1};\tx_{(i)i+1};\ldots;\tx_{(i)n}]\in\R^{n-1}.$$

We introduce the following estimation matrix:
$$
  \bx\triangleq\left(
     \begin{array}{ccc}
       \textrm{---}& x_{(1)}^\T & \textrm{---} \\
       \textrm{---}& x_{(2)}^\T & \textrm{---} \\
       &\vdots& \\
       \textrm{---}& x_{(n)}^\T & \textrm{---} \\
     \end{array}
   \right)\in\R^{n\times n},
$$
where $x^\T$ denotes the transpose of a column-vector $x$.
We let $\Om_a$ to denote an \an{augmented action set}, consisting of the vectors on the diagonals of the estimation matrices, i.e., $\Om_a = \{\bx\in\R^{n\times n}\,|\, \di(\bx) \in \Om\}$.
The pseudo-gradient estimation of the game is defined as $\bF(\bx)\in\mathbb R^{n\times n}$:
\begin{align}\label{eq:ext_gm}
   \bF(\bx)\triangleq\left(
     \begin{array}{cccc}
        \nabla_1 J_1(x_{(1)}) & 0 & \cdots &0\\
        0 & \nabla_2 J_2(x_{(2)}) & \cdots &0\\
        \vdots&\vdots&\ddots&\vdots\\
        0 & 0 & \cdots &\nabla_n J_n(x_{(n)})\\
     \end{array}
   \right).
\end{align}

The algorithm starts with an arbitrary initial $\bx^0\in\Om_a$, that is, each player $i$ holds an arbitrary point $x_{(i)}^0\in\R^{i-1}\times\Om_i\times\R^{n-i}$. All the subsequent estimation matrices $\bx_{1}, \bx_{2},\ldots$ are obtained
through the updates described by Algorithm 1.

\begin{remark}
    Algorithm~1 is inspired by the operator extrapolation approach presented  in~\cite{LanExtrapolation} for solving variational inequalities in a centralized setting. The extrapolation here corresponds to the expression $\nabla_i J_i(\hat{x}_{(i)}^k) + \lambda[\nabla_i J_i({x}_{(i)}^k) - \nabla_i J_i(\hat{x}_{(i)}^{k-1})]$ in the update of the individual actions $x_i^{k+1}$, $i\in[n]$. It is inspired by the connection between the Nesterov's acceleration and the gradients' extrapolation (see~\cite{Lan2018}). We call our proposed distributed algorithm \emph{Accelerated Direct Method} to emphasize that, in contrast to the work~\cite{AccGRANE_TAC}, the communication step $\hat{x}_{(i)}^k = \sum\limits_{j\in\N_{i}} w_{ij} x_{(j)}^k$ is directly implemented in the corresponding centralized procedure without any augmented game mapping.
\end{remark}

\begin{center}
  {\textbf{Algorithm 1: Accelerated Direct Method}}
  \smallskip

    \begin{tabular}{l}
    \hline
   \emph{  } Set mixing matrix $W$;\\
    \emph{  } Choose step size $\alpha>0$ and parameter $\lambda>0$;\\
    \emph{  } Pick arbitrary $x^0_{(i)}\in\R^{i-1}\times\Om_i\times\R^{n-i}$, $i=1,\ldots,n$; \\
    \emph{  } Set $\hat{x}_{(i)}^0 = \sum\limits_{j\in\N_i} w_{ij} x_{(j)}^0$ and $x_{(i)}^1 = \hat{x}_{(i)}^0$,\\
    \emph{  } \textbf{for} $k=1,2,\ldots$, all players $i=1,\ldots,n$ do \\
    \emph{  } \quad $\hat{x}_{(i)}^k = \sum\limits_{j\in\N_i} w_{ij} x_{(j)}^k$, \\
    \emph{  } \quad $x_{i}^{k+1}=\mathcal{P}_{\Om_i}
    \{\hat x_i^k-\alpha[\nabla_i J_i(\hat{x}_{(i)}^k)$\\
    \qquad\qquad\qquad\qquad\qquad\, $+ \lambda[\nabla_i J_i({x}_{(i)}^k) - \nabla_i J_i(\hat{x}_{(i)}^{k-1})]]\}$;\\
    \emph{  } \qquad\textbf{for} $\ell=\{1,\ldots,i-1,i+1,\ldots,n\}$ \\
    \emph{  } \qquad\qquad$\tx_{(i)\ell}^{k+1}=\hat{x}_{(i)\ell}^k$;\\
    \emph{  } \qquad\textbf{end for};\\
    \emph{~ } \textbf{end for}.\\
    \hline
    \end{tabular}
\end{center}

The compact expression of Algorithm~1 in terms of the pseudo-gradient estimation $\bF(\cdot)$ (see~\eqref{eq:ext_gm}) and the estimation matrices $\{\bx^k\}$ is as follows: 
\begin{align}\label{eq:compactAlg}
\hbx^k &= W\bx^k,\cr
\bx^{k+1} &= \pro{\Om_a}{\hbx^k - \alpha(\bF(\hbx^k)+\lambda(\bF(\bx^k) - \bF(\hbx^{k-1})))}.\qquad
\end{align}
In the further analysis of Algorithm~1, we will use its compact matrix form above. 

\subsection{Analysis}
\begin{lemma}\label{lem:Lip}
 Under Assumption~\ref{assum:Lipschitz}, the pseudo-gradient estimation $\bF(\cdot)$ is Lipschitz continuous over $\Om_a$ with the constant $L = \max_{i\in[n]}\sqrt{L^2_i + L^2_{-i}}$.
\end{lemma}
\proof
 See Lemma~1 in \cite{AccGRANE_TAC}. 
Before formulating the next result, let us define the following decomposition for each matrix $\bx\in\Om_a$: 
\an{
\begin{equation}\label{eq-perp-dec}
\bx = \bx_{||} + \bx_{\bot},\end{equation}
where $\bx_{||} = \frac{1}{n}\boldsymbol{1}\boldsymbol{1}^T\bx$ is the so called consensus matrix and $\bx_{\bot} = \bx - \bx_{||}$ with the implied property $\langle \bx_{||},\bx_{\bot} \rangle= 0$}.

Moreover, in the further analysis,  we will use the following lemma  (see Lemma 3.1 in \cite{LanBook} and the preceding discussion).
\begin{lemma}\label{lem:3p}
    Let $Y\subseteq \R^{n\times n}$ be a closed convex set, and let 
    $y_{t+1}\in Y$ be defined through the following relation: 
    \[y_{t+1}:=\mathcal{P}_Y [y_t-\alpha g_t]\]
    for some \an{$\alpha>0$ and  $g_t\in\R^{n\times n}$.}   Then, we have 
    \begin{align*}
        \gamma\langle g_t,& y_{t+1}-y\rangle + \frac{1}{2}\|y_t-y_{t+1}\|^2\cr&
        \le \frac{1}{2}\|y-y_{t}\|^2 - \frac{1}{2}\|y-y_{t+1}\|^2\quad\hbox{for all $y\in Y$}.
    \end{align*}
\end{lemma}


\begin{proposition}\label{prop:3points}
 Let Assumptions~\ref{assum:convex}-\ref{assum:connected} hold. Let $\bx^*$ be the consensus matrix with each row equal to the unique Nash equilibrium $x^*$ in the game $\Gamma$. Moreover, let a scalar $\eta\in(0,1)$ and a positive sequence $\{\theta_t\}$ be such that $\theta_{t-1}\eta(1-\eta)\ge \theta_tL^2\alpha^2\lambda^2$ for all $t\ge1$  and 
 $\theta_{t+1}\lambda = \theta_t$ for all $t\ge0$, 
 where $\lambda>0$ is the stepsize from Algorithm~1. Then, we have
 \begin{align}
 &\sum_{t=1}^k [\theta_t a_1 \left\|\bx^{t+1}_{||}  -\bx^*\right\|^2 + \theta_t a_2\|\bx^{t+1}_{\bot}\|^2]  \cr
  &\qquad- \frac{\theta_kL^2\alpha^2}{2(1-\eta)}\|\bx^{k+1} - \bx^*\|^2 \cr
  &\le \sum_{t=1}^k [\theta_t b_1 \left\|\bx^{t}_{||}  -\bx^*\right\|^2 + \theta_t b_2\|\bx^{t}_{\bot}\|^2], 
  \end{align}
where 
$a_1 = \frac{1-\eta}{2} + \frac{\mu}{2n}\alpha$, $a_2 = \frac{1-\eta}{2} - \left(L+\frac{2nL^2}{\mu}\right)\alpha$,	$b_1 =\frac{1+ \eta\|W-I\|^2}{2}$, $b_2 = \frac{(1-\eta)\sigma^2 + \eta(1+\|W-I\|^2)}{2}$. 
\end{proposition}
\proof
 \an{Applying Lemma~\ref{lem:3p}} to the iterates in \eqref{eq:compactAlg}, i.e., $y_t = \bx^t$, $Y = \Om_a$, \an{and 
 $g_t= \bF(\hbx^t)+\lambda(\bF(\bx^t) - \bF(\hbx^{t-1}))$,}
 we obtain
 $\alpha\left\langle \bF(\hbx^t)+\lambda(\bF(\bx^t) - \bF(\hbx^{t-1})), \bx^{t+1} - \bx^*\right\rangle + \frac{1}{2}\|\bx^{t+1} - \hbx^t\|^2\le\frac{1}{2}\|\hbx^t-\bx^*\|^2  - \frac{1}{2}\|\bx^{t+1} - \bx^*\|^2$ \an{for all $t\ge1$}, 
 \an{where $x^*$ is the Nash equilibrium from Assumption~\ref{assum:str_monotone}.}
  We multiply both sides of the preceding inequality by $\theta_t>0$, sum up the resulting relations over $t=1,\ldots, k,$ for an arbitrary $k\ge1$, to obtain
  {\allowdisplaybreaks
 \begin{align}\label{eq:3pl}
 &\sum_{t=1}^k[\theta_t\alpha\left\langle \bF(\hbx^t)+\lambda(\bF(\bx^t) - \bF(\hbx^{t-1})), \bx^{t+1} - \bx^*\right\rangle \cr
 &\qquad\qquad+ \frac{\theta_t}{2}\|\bx^{t+1} - \hbx^t\|^2]\cr
 &\le\sum_{t=1}^k\frac{\theta_t}{2}\left[\|\hbx^t-\bx^*\|^2  - \|\bx^{t+1} - \bx^*\|^2\right].
 \end{align}
  Next we consider the left hand side of the inequality above: 
 \begin{align}\label{eq:eq0}
  	&\sum_{t=1}^{k} \big[\theta_t \alpha\left\langle \bF(\hbx^t)+\lambda(\bF(\bx^t) - \bF(\hbx^{t-1})), \bx^{t+1} - \bx^*\right\rangle\cr
   &\qquad\qquad+ \frac{\theta_t}{2}\|\bx^{t+1} - \hbx^t\|^2\big]\cr
  	& = \sum_{t=1}^{k} 
	\big[\theta_t \alpha\langle \bF(\bx^{t+1}),\bx^{t+1} -\bx^*\rangle \cr
    &\quad -\theta_t\alpha \langle\bF(\bx^{t+1})-\bF(\hbx^{t}),\bx^{t+1} -\bx^*\rangle  \cr
    &\quad	+ \theta_t\alpha\lambda\langle\bF(\bx^{t})-\bF(\hbx^{t-1}),\bx^{t} -\bx^*\rangle \cr
  	&\quad+ \theta_t\alpha\lambda\langle\bF(\bx^{t})-\bF(\hbx^{t-1}),\bx^{t+1} -\bx^{t}\rangle + \frac{\theta_t}{2}\|\bx^{t+1} - \hbx^t\|^2\big]\cr
  	&=\sum_{t=1}^{k} 
	\big[\theta_t \alpha\langle\bF(\bx^{t+1}),\bx^{t+1} -\bx^*\rangle \cr
 &\quad+ \theta_t\alpha\lambda\langle\bF(\bx^{t})-\bF(\hbx^{t-1}),\bx^{t+1} -\bx^{t}\rangle + \frac{\theta_t}{2}\|\bx^{t+1} - \hbx^t\|^2\big]\cr
  	&\,-\theta_k\alpha \langle\bF(\bx^{k+1})-\bF(\hbx^{k}),\bx^{k+1} -\bx^*\rangle\cr
   &\,+\theta_0\alpha \langle\bF(\bx^{1})-\bF(\hbx^{0}),\bx^{1} -\bx^*\rangle\cr
   &=\sum_{t=1}^{k} 
	\big[\theta_t \alpha\langle\bF(\bx^{t+1}),\bx^{t+1} -\bx^*\rangle \cr
 &\quad+ \theta_t\alpha\lambda\langle\bF(\bx^{t})-\bF(\hbx^{t-1}),\bx^{t+1} -\bx^{t}\rangle + \frac{\theta_t}{2}\|\bx^{t+1} - \hbx^t\|^2\big]\cr
   &\,-\theta_k\alpha \langle\bF(\bx^{k+1})-\bF(\hbx^{k}),\bx^{k+1} -\bx^*\rangle, 
 \end{align}
 \an{where in the second equality we used
 $\theta_{1}\lambda=\theta_0$ and in the last equality we used $\bx^{1}=\hbx^{0}$}. 
\an{Next, we consider the sum of the second and the third term in relation~\eqref{eq:eq0}}
{\allowdisplaybreaks
\begin{align}\label{eq:eq1}
 &\sum_{t=1}^{k} \left[\theta_t\alpha\lambda\langle\bF(\bx^{t})-\bF(\hbx^{t-1}),\bx^{t+1} -\bx^{t}\rangle + \frac{\theta_t}{2}\|\bx^{t+1} - \hbx^t\|^2\right]\cr
 &=\sum_{t=1}^{k} \left[\theta_t\alpha\lambda\langle\bF(\bx^{t})-\bF(\hbx^{t-1}),\bx^{t+1} -\bx^{t}\rangle \right.\cr
 &\quad\left.+ \frac{\theta_t\eta}{2}\|\bx^{t+1} - \hbx^t\|^2 +\frac{\theta_{t-1}(1-\eta)}{2}\|\bx^{t} - \hbx^{t-1}\|^2\right]\cr
 &\qquad+ \frac{\theta_k(1-\eta)}{2}\|\bx^{k+1} - \hbx^k\|^2-\an{\frac{\theta_{0}(1-\eta)}{2}\|\bx^{1} - \hbx^{0}\|^2}\cr
 &=\sum_{t=1}^{k} \left[\theta_t\alpha\lambda\langle\bF(\bx^{t})-\bF(\hbx^{t-1}),\bx^{t+1} -\bx^{t}\rangle \right. \cr
 &\qquad\left.+\frac{\theta_t\eta}{2}\|\bx^{t+1} - \hbx^t\|^2+\frac{\theta_{t-1}(1-\eta)}{2}\|\bx^{t} - \hbx^{t-1}\|^2\right]\cr
 &\qquad+ \frac{\theta_k(1-\eta)}{2}\|\bx^{k+1} - \hbx^k\|^2\cr
 & \ge \sum_{t=1}^{k} 
 \left[-L\theta_t\alpha\lambda\|\bx^{t}-\hbx^{t-1}\|\|\bx^{t+1} -\bx^{t}\| + \frac{\theta_t\eta}{2}\|\bx^{t+1} - \bx^t\|^2 \right.\cr
 &\qquad\quad\left.+\frac{\theta_{t-1}(1-\eta)}{2}\|\bx^{t} - \hbx^{t-1}\|^2\right]\cr
 &\qquad + \sum_{t=1}^{k}\frac{\theta_t\eta}{2}
 \left[\|\bx^{t+1} - \hbx^t\|^2-\|\bx^{t+1} - \bx^t\|^2\right]\cr
 &\qquad+ \frac{\theta_k(1-\eta)}{2}\|\bx^{k+1} - \hbx^k\|^2\cr
 &\ge \sum_{t=1}^{k} 
 \left[-L\theta_t\alpha\lambda\|\bx^{t}-\hbx^{t-1}\|\|\bx^{t+1} -\bx^{t}\| + \frac{\theta_t\eta}{2}\|\bx^{t+1} - \bx^t\|^2 \right.\cr
 &\qquad\quad\left.+\frac{\theta_tL^2\alpha^2\lambda^2}{2\eta}\|\bx^{t} - \hbx^{t-1}\|^2\right]\cr
 &\qquad + \sum_{t=1}^{k}\frac{\theta_t\eta}{2}
 \left[\|\bx^{t+1} - \hbx^t\|^2-\|\bx^{t+1} - \bx^t\|^2\right]\cr
 &\qquad+ \frac{\theta_k(1-\eta)}{2}\|\bx^{k+1} - \hbx^k\|^2\cr
 &= \sum_{t=1}^{k} 
 (\|\bx^{t}-\hbx^{t-1}\|-\|\bx^{t+1} -\bx^{t}\|)^2\cr
 &\qquad + \sum_{t=1}^{k}\frac{\theta_t\eta}{2}
 \left[\|\bx^{t+1} - \hbx^t\|^2-\|\bx^{t+1} - \bx^t\|^2\right]\cr
 &\qquad+ \frac{\theta_k(1-\eta)}{2}\|\bx^{k+1} - \hbx^k\|^2\cr
 &\ge \sum_{t=1}^{k}\frac{\theta_t\eta}{2}[\|\bx^{t+1} - \hbx^t\|^2-\|\bx^{t+1} - \bx^t\|^2] \cr
 &\qquad+ \frac{\theta_k(1-\eta)}{2}\|\bx^{k+1} - \hbx^k\|^2,
\end{align}}
where in the first equality we used the relation $\bx^{1} = \hbx^{0}$}, in the first inequality we applied Lemma~\ref{lem:Lip} \an{and added and subtracted $(\theta_t\eta/2)\|\bx^{t+1} - \bx^t\|^2$}, whereas the second inequality is due to $\theta_{t-1}\eta(1-\eta)\ge \theta_tL^2\alpha^2\lambda^2$. 
Next, we use again Lemma~\ref{lem:Lip}  to obtain
{\allowdisplaybreaks
\begin{align}\label{eq:last}
 &-\theta_k\alpha \langle\bF(\bx^{k+1})-\bF(\hbx^{k}),\bx^{k+1} -\bx^*\rangle \cr
 &\qquad+ \frac{\theta_k(1-\eta)}{2}\|\bx^{k+1} - \hbx^k\|^2\cr
 &\ge - \theta_k\alpha L \|\bx^{k+1}-\hbx^{k}\|\|\bx^{k+1} -\bx^*\| \cr
 &\qquad+ \frac{\theta_k(1-\eta)}{2}\|\bx^{k+1} - \hbx^k\|^2\cr
 &\ge -\frac{L^2\theta_k\alpha^2}{2(1-\eta)}\|\bx^{k+1}-\bx^*\|^2.
\end{align}}
\tanya{The last inequality is obtained from $-2ab+a^2\ge-b^2$ with $a=\frac{\sqrt{\theta_k(1-\eta)}}{\sqrt{2}}\|\bx^{k+1} - \hbx^k\|$ and $b = \frac{L\sqrt{\theta_k}\alpha}{\sqrt{2(1-\eta)}}\|\bx^{k+1}-\bx^*\|$.}
The \an{preceding relations} imply
{\allowdisplaybreaks
\begin{align}\label{eq:eq2}
  	&\sum_{t=1}^{k} 
	\left[\theta_t \alpha\langle\bF(\bx^{t+1}),\bx^{t+1} -\bx^*\rangle + \frac{\theta_t}{2}\|\bx^{t+1} - \bx^*\|^2\right]\cr
 &\qquad-\frac{L^2\theta_k\alpha}{2(1-\eta)}\|\bx^{k+1}-\bx^*\|^2\cr
  &\stackrel{\eqref{eq:last}}{\le}\sum_{t=1}^{k} 
	\left[\theta_t \alpha\langle\bF(\bx^{t+1}),\bx^{t+1} -\bx^*\rangle + \frac{\theta_t}{2}\|\bx^{t+1} - \bx^*\|^2\right]\cr
 &\qquad-\theta_k\alpha \langle\bF(\bx^{k+1})-\bF(\hbx^{k}),\bx^{k+1} -\bx^*\rangle \cr
 &\qquad+ \frac{\theta_k(1-\eta)}{2}\|\bx^{k+1} - \hbx^k\|^2\cr
 &\stackrel{\eqref{eq:eq1}}{\le} \sum_{t=1}^{k} 
	\left[\theta_t \alpha\langle\bF(\bx^{t+1}),\bx^{t+1} -\bx^*\rangle + \frac{\theta_t}{2}\|\bx^{t+1} - \bx^*\|^2\right]\cr
 &\qquad-\theta_k\alpha \langle\bF(\bx^{k+1})-\bF(\hbx^{k}),\bx^{k+1} -\bx^*\rangle\cr
  &+ \sum_{t=1}^{k} \left[\theta_t\alpha\lambda\langle\bF(\bx^{t})-\bF(\hbx^{t-1}),\bx^{t+1} -\bx^{t}\rangle\right.\cr
 &\qquad\qquad\qquad\left. + \frac{\theta_t}{2}\|\bx^{t+1} - \hbx^t\|^2\right] \cr
  &\qquad - \sum_{t=1}^{k}\frac{\theta_t\eta}{2}[\|\bx^{t+1} - \hbx^t\|^2-\|\bx^{t+1} - \bx^t\|^2] \cr
  & \stackrel{\eqref{eq:eq0}}{=}\sum_{t=1}^{k} 
	\left[\frac{\theta_t}{2}\|\bx^{t+1} - \bx^*\|^2\right.\cr
 &\qquad\qquad\qquad\left.-\frac{\theta_t\eta}{2}[\|\bx^{t+1} - \hbx^t\|^2-\|\bx^{t+1} - \bx^t\|^2]\right]\cr
& \qquad + \sum_{t=1}^{k} \left[\theta_t \alpha\left\langle \bF(\hbx^t)+\lambda(\bF(\bx^t) - \bF(\hbx^{t-1})), \bx^{t+1} - \bx^*\right\rangle\right.\cr
&\qquad\qquad\qquad\left.+ \frac{\theta_t}{2}\|\bx^{t+1} - \hbx^t\|^2\right]\cr
   &\stackrel{\eqref{eq:3pl}}{\le} \sum_{t=1}^{k}
	\left[\frac{\theta_t}{2}\|\hbx^{t} - \bx^*\|^2 \right.\cr
 &\qquad\qquad\qquad\left.- \frac{\theta_t\eta}{2}\left[\|\bx^{t+1} - \hbx^t\|^2-\|\bx^{t+1} - \bx^t\|^2\right]\right]\cr
  	&=\sum_{t=1}^{k}
	\left[\frac{\theta_t(1-\eta)}{2}\|\hbx^{t} - \bx^*\|^2 \right.\\
 &\,\left.+ \frac{\theta_t\eta}{2}\left[\|\bx^{t+1} - \bx^t\|^2-\|\bx^{t+1} - \hbx^t\|^2+\|\hbx^{t} - \bx^*\|^2\right]\right].\nonumber
 \end{align}
}
\an{According to the definition~\eqref{eq:ext_gm} and taking into account that $x^*$ is the Nash equilibrium, we have }$\langle\bF(\bx^*),\bx-\bx^*\rangle = \langle F(x^*), x-x^*\rangle\ge 0$ for any $\bx\in\Om_a$ (\an{see~\eqref{eq:NE} and the definition of the mapping $\bF(\cdot)$ in~\eqref{eq:ext_gm}. 
Thus, we obtain
\begin{align*}
\langle\bF(\bx^{t+1}),\bx^{t+1} -\bx^*\rangle
&=\langle\bF(\bx^{t+1}) -\bF(\bx^*),\bx^{t+1} -\bx^*\rangle \cr
&\ +\langle\bF(\bx^*),\bx^{t+1}-\bx^*\rangle \cr
&\ge \langle\bF(\bx^{t+1})- \bF(\bx^*),\bx^{t+1} -\bx^*\rangle.
\end{align*}
}
\an{Recalling our notation for $\bx_{||}$ and $\bx_\perp$ (see~\eqref{eq-perp-dec}), } we have \tanya{$\|\bx^{t+1}  -\bx^*\| \le \|\bx_{||}^{t+1}  -\bx^*\| + \|\bx_{\bot}^{t+1}\|$}. 
\an{Next, by adding and substracting $\bF(\bx^{t+1}_{||})$ to the right hand side of the preceding inequality and using the relation $\bx^{t+1} = \bx_{||}^{t+1} +\bx_{\bot}^{t+1}$, we obtain 
\begin{align*}
&\langle\bF(\bx^{t+1}),\bx^{t+1} -\bx^*\rangle \ge\cr
 & \langle\bF(\bx^{t+1})-\bF(\bx^{t+1}_{||}) + \bF(\bx_{||}^{t+1}) - \bF(\bx^*),\cr
&\qquad\qquad\qquad\qquad\qquad\qquad\qquad
\bx_{||}^{t+1} +\bx_{\bot}^{t+1}  -\bx^*\rangle\cr
 &=\langle\bF(\bx^{t+1})- \bF(\bx_{||}^{t+1}),\bx^{t+1}  -\bx^*\rangle \cr
 &\qquad+ \langle
 \bF(\bx_{||}^{t+1})  - \bF(\bx^*),\bx_{||}^{t+1}  -\bx^*\rangle\cr
 &\qquad + \langle\bF(\bx_{||}^{t+1})  - \bF(\bx^*),\bx_{\bot}^{t+1}\rangle.\end{align*}
By Lemma~\ref{lem:Lip}, the mapping $\bF$ is Lipschitz continuous, implying that}
\begin{align}\label{eq:StrMon}
&\langle\bF(\bx^{t+1}),\bx^{t+1} -\bx^*\rangle \ge\cr
 &\ge -L\|\bx^{t+1}- \bx_{||}^{t+1}\|\|\bx^{t+1}  -\bx^*\|\cr
 &\qquad+\langle\bF(\bx_{||}^{t+1})  - \bF(\bx^*),\bx_{||}^{t+1}  -\bx^*\rangle \cr
 &\qquad- L\|\bx_{||}^{t+1}  - \bx^*\|\|\bx_{\bot}^{t+1}\|\cr
 & \ge -L\|\bx_{\bot}^{t+1}\|(\|\bx_{||}^{t+1}  -\bx^*\|+ \|\bx_{\bot}^{t+1}\|)\cr
 &\qquad + \langle \bF(\bx_{||}^{t+1})  - \bF(\bx^*),\bx_{||}^{t+1}  -\bx^*\rangle\cr
 &\qquad- L\|\bx_{||}^{t+1}  - \bx^*\|\|\bx_{\bot}^{t+1}\|\cr
 &\ge -2L\|\bx_{\bot}^{t+1}\|\|\bx_{||}^{t+1}  -\bx^*\| - L\|\bx_{\bot}^{t+1}\|^2 +\frac{\mu}{n}\|\bx_{||}^{t+1}  -\bx^*\|^2 \cr
 &\ge \frac{\mu}{2n}\|\bx_{||}^{t+1}  -\bx^*\|^2 - \left(L+\frac{2nL^2}{\mu}\right)\|\bx_{\bot}^{t+1}\|^2,
\end{align}
where in the last two inequalities we used Lemma~\ref{lem:Lip} and Assumption~\ref{assum:str_monotone}, implying that $\langle\bF(\bx_{||}^{t+1})  - \bF(\bx^*),\bx_{||}^{t+1}  -\bx^*\rangle 
= \langle F(x_{||}^{t+1}) - F(x^*), x_{||}^{t+1}  -x^*\rangle
\ge \mu\|x_{||}^{t+1}  -x^*\|^2 
= \frac{\mu}{n}\|\bx_{||}^{t+1}  -\bx^*\|^2$,
and the fact that $-2L\|\bx_{\bot}^{t+1}\|\|\bx_{||}^{t+1}  -\bx^*\| \ge -\frac{2nL^2}{\mu}\|\bx_{\bot}^{t+1}\|^2 - \frac{\mu}{2n}\|\|\bx_{||}^{t+1}  -\bx^*\|^2$. 
As for the right hand side of~\eqref{eq:eq2}, we notice that 
\begin{align}\label{eq:rhs1}
 &\|\hbx^{t} - \bx^*\|^2 
 = \|W(\bx_{||}^{t} + \bx_{\bot}^{t}) - \bx^*\|^2 \cr
 &= \|\bx_{||}^{t}- \bx^*\|^2 + \|W\bx_{\bot}^t\|^2 \cr
 &\le \|\bx_{||}^{t}- \bx^*\|^2 + \sigma^2\|\bx_{\bot}^t\|^2,
\end{align}
where in the last inequality we used~\eqref{eq:sigma}.
Moreover,
\begin{align}\label{eq:rhs2}
 &\|\bx^{t+1} - \bx^t\|^2-\|\bx^{t+1} - \hbx^t\|^2+\|\hbx^{t} - \bx^*\|^2 \cr
 &\le \|\bx^t - \bx^*\|^2 + \|\hbx^t-\bx^t\|^2 + \|\bx^* - \bx^{t+1}\|^2.
\end{align}
Combining~\eqref{eq:StrMon}-\eqref{eq:rhs2} with~\eqref{eq:eq2} and taking into account that 
\[\|\hbx^t-\bx^t\|^2 
= \|W\bx^t-\bx^t-W\bx^*+\bx^*\|^2 
= \|(I-W)(\bx^t-\bx^*)\|,\] 
we obtain 
\begin{align}\label{eq:final}
 &\sum_{t=1}^{k} 
 \left[\theta_t \alpha\left(\frac{\mu}{2n}\|\bx_{||}^{t+1}  -\bx^*\|^2 - \left(L+\frac{2nL^2}{\mu}\right)\|\bx_{\bot}^{t+1}\|^2\right)\right. \cr
 &\quad\left.+ \frac{\theta_t(1-\eta)}{2}\|\bx^{t+1} - \bx^*\|^2\right]-\frac{L^2\theta_k\alpha^2}{2(1-\eta)}\|\bx^{k+1}-\bx^*\|^2\cr
 &\le \sum_{t=1}^{k}\frac{\theta_t(1-\eta)}{2}
 \left[\|\bx_{||}^{t}- \bx^*\|^2 + \sigma^2\|\bx_{\bot}^t\|^2\right] \cr
 &\,+\sum_{t=1}^{k}\frac{\theta_t\eta}{2}
 \left[\|\bx_{||}^{t}- \bx^*\|^2 + \|\bx_{\bot}^t\|^2\right](1+\|(I-W)\|^2).
\end{align}
After rearranging the terms we conclude the result. 

Next we formulate our main result. 

\begin{theorem}\label{th:main}
	Let the parameters in the Algorithm~1 be chosen as follows: 
	\begin{align}\label{eq:g} 
		&\alpha\le \min\{g_1, g_2, g_3, g_4\}, \quad \lambda = \frac{1}{1+\epsilon(\alpha)},
	\end{align}
where 
\begin{align*}
    &g_1 =  \frac{n\mu(1-\sigma^2)}{4(\mu+2nL)^2(1+\dW)},\cr
	&g_2 = \frac{n(1+\dW)}{2\mu},\cr
	&g_3 = \frac{\mu n(1+\dW)}{\mu^2+L^2(1+\dW)^2n^2}, \cr
	&g_4 = \frac{\mu}{\sqrt{4L^2\mu^2+16(L\mu+2nL^2)^2}},\cr
	&\epsilon(\alpha) = \frac{2\mu\alpha/n - (1+\dW)(1-\sqrt{1-4L^2\alpha^2})}{2+\dW(1-\sqrt{1-4L^2\alpha^2})}.	
\end{align*}
Then $\epsilon(\alpha)>0$ and 
\begin{align*}
   &\|\bx^{k+1} - \bx^*\|^2\cr
   &\le \frac{8+4\dW-4\dW\sqrt{1-4L^2\alpha^2}}{(1+\epsilon(\alpha))^{k-1}}\|\bx^1-\bx^*\|^2.
\end{align*}
\end{theorem}
\proof
  	Let $\theta_t=c^t$, where 
  	\begin{align}\label{eq:c}
  		c = \min\left\{\frac{a_1}{b_1},\frac{a_2}{b_2}\right\},
  	\end{align}
 and $a_1$, $a_2$, $b_1$, and $b_2$ are defined in Proposition~\ref{prop:3points}, namely,
 \begin{align*}
 	&a_1 = \frac{1-\eta}{2} + \frac{\mu}{2n}\alpha, \quad a_2 = \frac{1-\eta}{2} - \left(L+\frac{2nL^2}{\mu}\right)\alpha, \cr
 	& b_1 =\frac{1+ \eta\|W-I\|^2}{2} , \quad b_2 = \frac{(1-\eta)\sigma^2 + \eta(1+\|W-I\|^2)}{2}. 
 \end{align*}
Moreover, let us choose $\lambda = \frac{\theta_t}{\theta_{t+1}}$ and $\eta = \frac{1-\sqrt{1-4L^2\alpha^2}}{2}$. 
Next, we demonstrate that under the condition~\eqref{eq:g}, $c>1$. For this purpose we check that in this case both $\frac{a_1}{b_1}$ and $\frac{a_2}{b_2}$ are larger than 1. 
Indeed, 
\begin{align*}
	\frac{a_1}{b_1} &= \frac{1-\eta+\mu\alpha/n}{1+\eta\|W-I\|^2},\cr
    \frac{a_2}{b_2} &= \frac{1-\eta-2\left(L+\frac{2nL^2}{\mu}\right)\alpha}{(1-\eta)\sigma^2 + \eta(1+\|W-I\|^2)}.
\end{align*}
First, we notice that under the condition $\alpha\le g_1$,  
$\frac{a_1}{b_1}\le\frac{a_2}{b_2}$ (see Appendix). Next,
\[c=\frac{a_1}{b_1}>1,\]
 if and only if 
 \begin{align}\label{eq:eta}
 	\eta&<\frac{\mu\alpha}{n(1+\|W-I\|^2)}.
 \end{align}
Given the definition of $\eta$, we conclude that 
under the condition $\alpha< \min\{g_2,g_3\}$,
\[\eta = \frac{1-\sqrt{1-4L^2\alpha^2}}{2} < \frac{\mu\alpha}{n(1+\|W-I\|^2)},\]
and, thus, \eqref{eq:eta} holds. 

\an{Since $c>1$, it follows that}
\begin{align*}
	L^2\alpha^2 = \eta(1-\eta)\le c \eta(1-\eta), 
\end{align*}
which implies $\theta_{t-1}\eta(1-\eta)\ge \theta_tL^2\alpha^2\lambda^2$. Hence, the conditions of Proposition~\ref{prop:3points} hold and we conclude that 
 \begin{align*}
	&\theta_k[a_1 \|\bx_{||}^{k+1}  -\bx^*\|^2 +  a_2\|\bx_{\bot}^{t+1}\|^2]  - \frac{\theta_kL^2\alpha^2}{2(1-\eta)}\|\bx^{k+1} - \bx^*\|^2 \cr
 &\le \theta_1 [b_1 \|\bx_{||}^{1}  -\bx^*\|^2 +  b_2\|\bx_{\bot}^{1}\|^2], 
\end{align*}
where we used definition of $c$ in~\eqref{eq:c} implying $\theta_t b_1\le \theta_{t-1} a_1$ and $\theta_t b_2\le \theta_{t-1} a_2$. \an{We also use the fact that $\bx^*_\bot=0$.}
Thus, 
 $	\theta_k\min\{a_1,a_2\}\|\bx^{k+1}  -\bx^*\|^2  - \frac{\theta_kL^2\alpha^2}{2(1-\eta)}\|\bx^{k+1} - \bx^*\|^2 
	\le \theta_1 \max\{b_1,b_2\} \|\bx^{1}  -\bx^*\|^2.
$
As $\min\{a_1,a_2\} = a_2$ and $\max\{b_1,b_2\} = b_1$, 
we conclude that 
 \begin{align}\label{eq:main1}
	\theta_k\left[a_2- \frac{L^2\alpha^2}{2(1-\eta)}\right]\|\bx^{k+1}  -\bx^*\|^2  \le \theta_1 b_1 \|\bx^{1}  -\bx^*\|^2.
	\end{align}
Next, we notice that $a_2 = \frac{1-\eta}{2} - \left(L+\frac{2nL^2}{\mu}\right)\alpha$ is larger or equal to $1/8$, if the conditions $\alpha\le\frac{\mu}{4(L\mu+2nL^2)}$ and $\alpha\le\frac{\sqrt 3}{4L}$ hold (see Appendix).
On the other hand, given the condition $\alpha \le \frac{\sqrt 7}{8L}$, we have $\frac{L^2\alpha^2}{2(1-\eta)}\le \frac{1}{16}$. Taking this inequality together with $a_2\ge \frac{1}{8}$ into account, we obtain from~\eqref{eq:main1} that, given $\alpha\le\min\{g_4,g_5\}$, 
 \begin{align}\label{eq:main2}
	\frac{\theta_k}{16}\|\bx_{k+1}  -\bx^*\|^2  \le \theta_1 b_1 \|\bx_{1}  -\bx^*\|^2.
\end{align}
Finally, we use that $\theta_k = c^k = (1+\epsilon(\alpha))^k$, where 
$\epsilon(\alpha) = \frac{a_1}{b_1} - 1 = \frac{2\mu\alpha/n - (1+\dW)(1-\sqrt{1-4L^2\alpha^2})}{2+\dW(1-\sqrt{1-4L^2\alpha^2})}$
to conclude the result from~\eqref{eq:main2}. 

\begin{corollary}\label{cor:results}
	Taking into account the conditions of the theorem above, one can choose the step size $\alpha = O\left(\frac{\mu}{L^2n}\right)$ to obtain the convergence rate $O\left(\exp\left\{-\frac{k}{\gamma^2n^2}\right\}\right)$ which is faster than the rates of previously proposed methods for distributed learning of Nash equilibria in restricted strongly monotone games~\cite{AccGRANE_TAC, Bianchi2019, ifac_TatNed}. Indeed, the GRANE algorithm from~\cite{AccGRANE_TAC} is proven to converge to the Nash equilibrium with the rate $O\left(\exp\left\{-\frac{k}{\gamma^6n^6}\right\}\right)$, whereas the direct distributed procedure in~\cite{Bianchi2019, ifac_TatNed} improves this rate to $O\left(\exp\left\{-\frac{k}{\gamma^4n^3}\right\}\right)$.
\end{corollary}

\section{Simulations}
\begin{figure}[!htb]
	\centering
	\includegraphics[width=0.5\textwidth]{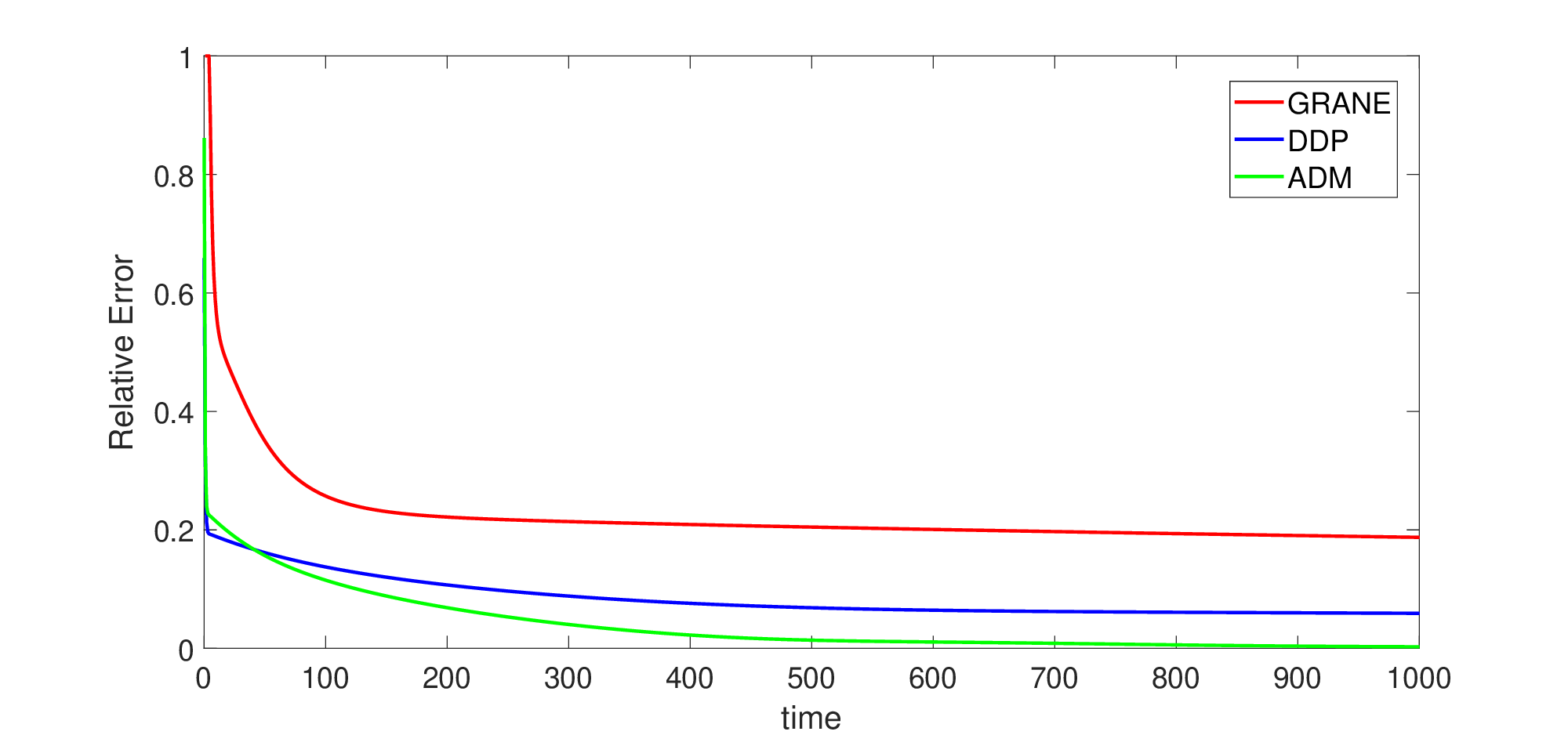}
	\caption{\label{eps:compare}Comparison of the presented accelerated direct method (ADM) with GRANE and DDP}
\end{figure}
Let us consider a class of games  with strongly monotone game mappings. Specifically, we have players $\{1,2,\ldots,20\}$ and each player $i$'s objective is to minimize the cost function $J_i(x_i,x_{-i})=f_i(x_i)+l_i(x_{-i})x_i$, where $f_i(x_i)=0.5a_ix_i^2+b_ix_i$ and $l_i(x_{-i})=\sum_{j\neq i}c_{ij}x_j$. The local cost function is dependent on actions of all players, but the underlying communication graph is a randomly generated tree graph. We randomly select $a_i>0$, $b_i$, and $c_{ij}$ for all possible $i$ and $j$ to guarantee strong monotonicity of the pseudo-gradient.

We simulate the proposed gradient play algorithm and compare its implementation with the implementations of the algorithm GRANE presented in \cite{AccGRANE_TAC} and direct distributed procedure (DDP) from \cite{Bianchi2019, ifac_TatNed}. Figure~\ref{eps:compare} demonstrates the simulation results which support theoretic ones stated in Corollary~\ref{cor:results}.

\section{Conclusion}
This work extends centralized operator extrapolation method presented in \cite{LanExtrapolation} to distributed settings in restricted strongly monotone games where players can exchange their information only with local neighbors via some communication graph. The proposed procedure is proven to possess a geometric rate and to outperform the previously developed algorithms calculating Nash equilibria in games under the same assumptions. Future research directions include consideration of a more general communication topology and study of \an{lower} bounds for convergence rates of distributed methods in such class of games.

\bibliographystyle{plain}
\bibliography{accNE}

\appendix
\emph{1. More details on the inequality $\frac{a_1}{b_1}\le\frac{a_2}{b_2}$ under the condition $\alpha<g_1$.}
{\allowdisplaybreaks
The condition $\alpha<g_1 = \frac{n\mu(1-\sigma^2)}{4(\mu+2nL)^2(1+\dW)}$ implies
\begin{align*}
 \alpha	\left(1+\frac{2nL}{\mu}\right)^2 \le \frac{n(1-\sigma^2)}{4\mu(1+\dW)}.
\end{align*}
Thus, since $\eta = \frac{1-\sqrt{1-4L^2\alpha^2}}{2}<\frac{1}{2}$, we conclude that $(1-\eta)^2>\frac{1}{4}$, and, hence
\begin{align*}
	\alpha	\left(1+\frac{2nL}{\mu}\right)^2(1+\dW) \le \frac{n}{\mu}(1-\sigma^2)(1-\eta)^2.
\end{align*}
Next, we use $\eta<1$ and $\sigma<1$ to obtain
\begin{align*}&\left(1+\frac{2nL}{\mu}\right)^2(1+\dW) \cr
&\ge \left(1-\frac{2nL}{\mu}\right)^2(1+\eta\dW) - (1-\eta)(1-\sigma^2).
\end{align*}
Combining two last inequalities, we get 
\begin{align*}
	&\alpha	[\left(1+\frac{2nL}{\mu}\right)^2(1+\eta\dW)- (1-\eta)(1-\sigma^2)] \cr
 &\le \frac{n}{\mu}(1-\sigma^2)(1-\eta)^2.
\end{align*}
By multiplying both sides by $\frac{\mu}{n}$, we obtain
\begin{align*}
	&\alpha	\left[\left(\frac{\mu}{n}+2L+\frac{4nL^2}{\mu}\right)(1+\eta\dW)\right.\cr
 &\qquad\left.- (1-\eta)(1-\sigma^2)\frac{\mu}{n}\right] \le (1-\sigma^2)(1-\eta)^2\cr
	&\qquad\qquad\Updownarrow\cr
	&(1-\eta)(1+\eta\dW)+\alpha\frac{\mu}{n}(1+\eta\dW) \cr
 &\qquad\qquad- (1-\eta)^2(1-\sigma^2) - \frac{\mu}{n}\alpha(1-\eta)(1-\sigma^2)\cr
	&\le
	(1-\eta)(1+\eta\dW)\cr
 &\qquad\qquad-(1+\eta\dW)2\left(L+\frac{2nL^2}{\mu}\right)\alpha\cr
	&\qquad\qquad\Updownarrow\cr
	&(1+\eta\dW)\left(1-\eta+\alpha\frac{\mu}{n}\right)\cr
 &\qquad\qquad-(1-\eta)^2(1-\sigma^2)\left(1-\eta+\alpha\frac{\mu}{n}\right)\cr
	&\le(1+\eta\dW)\left(1-\eta-2\left(L+\frac{2nL^2}{\mu}\right)\alpha\right)\cr
	&\qquad\qquad\Updownarrow\cr
	&(1+\eta\dW-(1-\eta)^2(1-\sigma^2))\left(1-\eta+\alpha\frac{\mu}{n}\right)\cr
 &\le(1+\eta\dW)\left(1-\eta-2\left(L+\frac{2nL^2}{\mu}\right)\alpha\right)\cr
		&\qquad\qquad\Updownarrow\cr
	&\frac{1-\eta-2\left(L+\frac{2nL^2}{\mu}\right)\alpha}{(1-\eta)\sigma^2 + \eta(1+\|W-I\|^2)}\le\frac{1-\eta+\mu\alpha/n}{1+\eta\|W-I\|^2}.
\end{align*}}
\medskip

\emph{2. More details on the inequality  $a_2 = \frac{1-\eta}{2} - \left(L+\frac{2nL^2}{\mu}\right)\alpha\ge1/8$ under the conditions $\alpha\le\frac{\mu}{4(L\mu+2nL^2)}$ and $\alpha\le\frac{\sqrt 3}{4L}$.}

As $\eta = \frac{1-\sqrt{1-4L^2\alpha^2}}{2}$,
\[\frac{1-\eta}{2} - \left(L+\frac{2nL^2}{\mu}\right)\alpha\ge1/8\]
if and only if 
\begin{align*}
	&2+2\sqrt{1-4L^2\alpha^2} - 8\left(L+\frac{2nL^2}{\mu}\right)\alpha\ge1\cr
	&\qquad\qquad\Updownarrow\cr
	&2 - 8\left(L+\frac{2nL^2}{\mu}\right)\alpha\ge1-2\sqrt{1-4L^2\alpha^2},
\end{align*}
which holds, if $2 - 8\left(L+\frac{2nL^2}{\mu}\right)\alpha\ge0$ and  $1-2\sqrt{1-4L^2\alpha^2}\le0$. 
The first inequality is guaranteed by $\alpha\le\frac{\mu}{4(L\mu+2nL^2)}$, whereas the second one is implied by $\alpha\le\frac{\sqrt 3}{4L}$.

\end{document}